\documentclass[12pt]{article}
\usepackage{graphicx} % Required for inserting images

\usepackage{geometry}
\usepackage{parskip}

\usepackage{amsmath}
\usepackage{amssymb}

\usepackage{mathtools}

\title{Some Identities For Periods of Hulek-Verrill Threefolds}
\author{Xenia de la Ossa, Mohamed Elmi}
\date{}

\begin{document}
\maketitle

We study the Hulek-Verrill families of Calabi-Yau threefolds. They are birationally equivalent to fibred products of elliptic surfaces, so we expect to be able to compute periods on these threefolds by integrating products of elliptic periods over a contour on $\mathbb{P}^1$. We numerically verify this in several examples. This article was submitted to MATRIX Annals (2024) for inclusion in the proceedings of the conference “The Geometry of Moduli Spaces in String Theory”, held 2–13 September 2024. %The Hulek-Verril threefolds are interesting because some of them are attractor varieties of rank two. These are modular threefolds where the Hodge structure splits over the rational numbers in a specific way. This note is part of an ongoing effort to better understand attractor varieties of rank two.

% The Hulek-Verril threefolds are birationally equivalent to families of 

% In this note, we begin to take steps towards extending the work by one of authors \cite{Elmi:2023hof} to Calabi-Yau threefolds which are not necessarily Hadamard products, but are also related to the fibred product of elliptic surfaces. In particular, we work with the Hullek-Verrill example of \cite{2003math......4169H} and \cite{2004math......7327V}. In this family of Calabi-Yau threefolds, there are a number of attractor points of rank two where the Hodge structure splits. The ultimate goal would be to understand D-branes and special cycles at these points. In this note, we compute periods of the threefolds as integrals of products of elliptic periods.

%tableofcontents
\section{Hulek-Verril family of $1$-folds $\mathcal{E}_{1,1,\frac{1}{\varphi}}$}

The Hulek-Verril elliptic surface $\mathcal{E}_{a,b,c}$ is the hypersurface in $\mathbb{P}^2\times\mathbb{P}^1$, defined by
\begin{equation}
    P(X,\lambda)=(X^1+X^2+X^3)(a X^2X^3+b X^1X^3+ c X^1X^2)\lambda^0-\lambda^1X^1X^2X^3=0
\end{equation}
where $[X^1,X^2,X^3]\in\mathbb{P}^2$, $[\lambda^0,\lambda^1]=\mathbb{P}^1$ and $(a,b,c)\in\mathbb{C}^3$ are parameters.\cite{2003math......4169H}

We are interested in the elliptic surface with parameters $\mathcal{E}_{1,1,\frac{1}{\varphi}}$. We choose the affine coordinate $\lambda=\frac{\lambda^0}{\lambda^1}$ on $\mathbb{P}^1$ and work on the coordinate patch $X^3=1$ on $\mathbb{P}^2$. With these choices, the holomorphic $1$-form on the elliptic fibre at $\lambda$ is given by
\begin{equation}
    \frac{1}{2\pi i}\oint_{C} \frac{dX^1\wedge dX^2}{P(X,\lambda)}
\end{equation}
where $C$ is a ``small contour" that encircles the $P(X,\lambda)=0$ locus. 

The holomorphic period at $\lambda=0$ is given by
\begin{equation}
     \frac{1}{2\pi i}\oint_{|X^1|=\epsilon}\oint_{|X^2|=\delta} \frac{dX^1\wedge dX^2}{P(X,\lambda)}=(2\pi i)f_0(\lambda,\varphi)
\end{equation}
for sufficiently small $\epsilon$ and $\delta$. The function $f_0(\lambda,\varphi)$ given by
\begin{align}
\begin{split}
\label{eq: f_0 for HV elliptic curve}
    f_0(\lambda,\varphi)&=\sum_{n=0}^\infty\left(
    \sum_{k=0}^n\binom{n}{k}^2\binom{2k}{k}\varphi^{-n+k}\right)\lambda^n\\
    % &=1+\left(2+\frac{1}{\varphi}\right)\lambda+\left(6+\frac{8}{\varphi}+\frac{1}{\varphi^2}\right)\lambda^2+\left(20+\frac{54}{\varphi}+\frac{18}{\varphi^2}+\frac{1}{\varphi^3}\right)\lambda^3+\mathcal{O}(\lambda^4).
\end{split}
\end{align}
%\begin{equation}
%\label{eq: f_0 for HV elliptic curve}
%    f_0(\lambda,\varphi)= 1+\left(2+\frac{1}{\varphi}\right)\lambda+\left(6+\frac{8}{\varphi}+\frac{1}{\varphi^2}\right)\lambda^2+\left(20+\frac{54}{\varphi}+\frac{18}{\varphi^2}+\frac{1}{\varphi^3}\right)\lambda^3+\mathcal{O}(\lambda^4).
%\end{equation}
and the coefficient of $\lambda^n$ is given by the constant term of
\begin{equation}
    \left[ (X^1+X^2+1)\left(\frac{1}{X^1}+\frac{1}{X^2}+\frac{1}{\varphi}\right)\right]^n.
\end{equation}

Both periods of the elliptic curve at $\lambda$ are annihilated by the Picard-Fuchs differential operator
\begin{equation}
\label{eq: second order PF operator at generic varphi}
    \mathcal{L}=R_2(\lambda,\varphi)(\theta+1)^2+R_1(\lambda,\varphi)(\theta+1)+R_0(\lambda,\varphi)
\end{equation}
where $\theta=\lambda\frac{d}{d\lambda}$ and 
\begin{align}
\begin{split}
    R_2(\lambda,\varphi)&=\left(3+\lambda\left(\frac{1}{\varphi}-4\right)\right)\left(1-\frac{\lambda}{\varphi}\right)\left(1-2\lambda\left(\frac{1}{\varphi}+4\right)+\lambda^2\left(\frac{1}{\varphi}-4\right)^2\right)\\
    R_1(\lambda,\varphi)&=-6+6\lambda\left(\frac{1}{\varphi}+6\right)+2\lambda^2\left(\frac{1}{\varphi}-4\right)\left(\frac{3}{\varphi}+8\right)-2\lambda^3\left(\frac{1}{\varphi}-4\right)\left(8-\frac{6}{\varphi}+\frac{3}{\varphi^2}\right)\\
    R_0(\lambda,\varphi)&=3-\lambda\left(\frac{1}{\varphi}+14\right)-\lambda^2\left(\frac{1}{\varphi}+2\right)\left(\frac{3}{\varphi}-10\right)+\lambda^3\left(\frac{1}{\varphi}-4\right)\left(4-\frac{2}{\varphi}+\frac{1}{\varphi^2}\right)
\end{split}
\end{align}

This differential operator $\mathcal{L}$ is implicitly given by the recursion relations of Verril in \cite{2004math......7327V} and it has the Riemann symbol
\begin{equation}
    \mathcal{P}~\left\{~~\begin{array}{c c c c c} 
    ~\frac{3\varphi}{4\varphi-1}~ & ~0 ~& ~\varphi~ & ~\frac{\varphi}{(1\pm2\sqrt{\varphi})^2}~ & ~\infty~
    \\[8pt]
    \hline
    \\[-5pt]
    0 & 0 & 0 & 0 & 1\\
    2 & 0 & 0 & 0 & 1\\
    \end{array}~~; ~\lambda~\right\}~.
\end{equation}

For generic values of $\varphi$, the elliptic periods $\omega(s,\varphi)$ are annhilated by the Picard-Fuchs operator \eqref{eq: second order PF operator at generic varphi}. However, at $\varphi=1$, some of the singularities of \eqref{eq: second order PF operator at generic varphi} coalesce and $\omega(\lambda,1)$ is annihilated by a Picard-Fuchs operator with Riemann symbol
\begin{equation}
    \mathcal{P}~\left\{~\begin{array}{c c c c}  0 & \frac{1}{9} & 1 & \infty\\[5pt]
    \hline
    \\[-5pt]
     0 & 0 & 0 & 1\\
     0 & 0 & 0 & 1\\
    \end{array}~~; ~\lambda\right\}~.
\end{equation}
This will be important in the following section. We will will also need a basis of periods with integral symplectic monodromy. After some experimentation, we make the following conjecture.

\textbf{Conjecture:} The basis of periods
\begin{equation}
\label{eq: integral symplectic basis of HV ellipitc curves}
    \omega(\lambda,\varphi)
    =(2\pi i)
    \begin{pmatrix} 1 & 0 \\[3pt]
    -\frac{\log (\varphi)}{2 \pi i} & \frac{3}{2\pi i} \end{pmatrix}
    \left(\begin{array}{l}
    f_0(\lambda,\varphi) \\[3pt] 
    f_0(\lambda,\varphi) \log(\lambda) +f_1(\lambda,\varphi)
    \end{array}\right)
\end{equation}
has monodromy in $SL_2(\mathbb{Z})$ for all $\varphi\in\mathbb{C}^*$, where $f_0(\lambda,\varphi)$ is given in \eqref{eq: f_0 for HV elliptic curve} and $f_1(\lambda,\varphi)$ is computed from the Picard-Fuchs equation \eqref{eq: second order PF operator at generic varphi} and is normalised so that $f_1(0,\varphi)=0$. \footnote{We are grateful to Joseph McGovern for helpful discussions related to this choice of basis.}

\section{Hulek-Verril family of threefolds as a fibred product of elliptic surfaces}
AESZ 34 describes the variation of the Hodge structure of (a quotient of) the one-parameter fibred products of elliptic surfaces
\begin{align}
\label{eq: one parameter HV threefold}
    \begin{split}
        P_1(X,\lambda)&=(X^1 + X^2 + X^3)(X^2X^3 + X^1X^3 + X^1X^2)\lambda^0-\lambda^1 X^1X^2X^3 =0\\
        P_2(Y,\lambda) &=(Y^1 + Y^2 + Y^3)\left(Y^2Y^3 + Y^1Y^3 + \frac{1}{\varphi} Y^1Y^2\right)\lambda^0-\lambda^1 Y^1Y^2Y^3=0
    \end{split}
\end{align}
where $\varphi\in\mathbb{C}$ is a parameter. These equations describe the fibred product of Hulek-Verril elliptic surfaces $\mathcal{E}_{1,1,1}\times_{\mathbb{P}^1}\mathcal{E}_{1,1,\frac{1}{\varphi}}$. Equivalently, these equations describe a one-parameter family of complete intersections with configuration matrix
\begin{equation}
\begin{matrix}
    \mathbb{P}^2\\
\mathbb{P}^2\\
\mathbb{P}^1
\end{matrix}
\begin{bmatrix}
    3 & 0 \\
   0 & 3 \\
   1 & 1
\end{bmatrix}\, .
\end{equation}
On the $X^3=Y^3=\lambda^1=1$ coordinate patch, the holomorphic $3$-form takes the form
\begin{equation}
    \Omega = \frac{1}{(2\pi i)^2}\int_C \frac{dX_1\wedge dX_2\wedge dY_1\wedge dY_2\wedge d\lambda_0}{P_1(X,\lambda)P_2(Y,\lambda)}\, ,
\end{equation}
where $C$ is a ``small contour" that encircles the $P_1(X,\lambda)=P_2(Y,\lambda)=0$ locus.

Let $\lambda=\frac{\lambda^0}{\lambda^1}$ and $G_i\in H_1(E_1,\mathbb{Z})\otimes H_1(E_2,\mathbb{Z})$ where $E_1$ and $E_2$ are the elliptic fibres defined by $P_1(X,\lambda)=0$ and $P_2(Y,\lambda)=0$ respectively. For suitably chosen $G_i$ and contours $\ell_{i}$ in the complex $\lambda$ plane, we expect an identify
\begin{equation}
\label{eq: generic integral of elliptic periods}
    \Gamma^T\Sigma_4\Pi(\varphi)=\sum_{i=1}^N G_i^T\Sigma_{2,2}\int_{\ell_i} \omega(\lambda,1)\otimes \omega(\lambda,\varphi)~d\lambda\, ,
\end{equation}
for some $\Gamma\in\mathbb{Q}^4$ and $\Pi$ that is, up a gauge transformation, a solution of AESZ 34 with integral symplectic monodromy \cite{2005math......7430A}.  We have identified $G_i$ with a vector of periods in $\mathbb{Z}^4$ and defined
\begin{equation} \Sigma_n=\begin{pmatrix} 0 & I_{\frac{n}{2}}\\ -I_{\frac{n}{2}} & 0 \end{pmatrix}
\end{equation}
and
\begin{equation}
    \Sigma_{m,n}=\Sigma_m\otimes\Sigma_n.
\end{equation}

A basis of solutions of AESZ 34 with integral symplectic monodromy is given by
\begin{equation}
\label{eq: change of basis matrix at MUM point}
    \Pi(\varphi)
    = (2\pi i)^3\begin{pmatrix} 
    ~\frac{\zeta(3)}{(2\pi i)^3}\,\chi ~& \frac{c_2\cdot H}{24} & ~0~ & \frac{H^3}{3!}\\[3pt]
    \frac{c_2\cdot H}{24} & \frac{\sigma}{2} & -\frac{H^3}{2!} & 0\\[3pt]
    1 & 0 & 0 & 0\\
    0 & 1 & 0 & 0
    \end{pmatrix}\begin{pmatrix} 1 & 0 & 0 & 0\\ 0 & (2\pi i)^{-1} & 0 & 0\\ 0 & 0 & (2\pi i)^{-2} & 0\\ 0 & 0 & 0 & (2\pi i)^{-3}\end{pmatrix}\varpi(\varphi)\, ,
\end{equation}
where $(H^3,c_2\cdot H,\chi)=(12,12,-8)$, $\sigma=H^3~\text{mod} ~2=0$ and $\varpi=\left(\varpi_0,\varpi_1,\varpi_2,\varpi_3\right)$. The components of $\varpi$ are 
% \begin{equation}
% \label{eq: Frobenius solutions at varphi=0}
%     \varpi(\varphi)=\varphi\left(\begin{array}{l}
%     F_0(\varphi) \\ 
%     F_0(\varphi) \log(\varphi) +F_1(\varphi)\\ 
%     F_0(\varphi) \log^2(\varphi) +2F_1(\varphi)\log(\varphi)+F_2(\varphi)\\
%     F_0(\varphi)\log^3(\varphi)+3F_1(\varphi)\log^2(\varphi)+3F_2(\varphi)\log(\varphi)+F_3(\varphi)
%     \end{array}\right)
% \end{equation}
\begin{equation}
\label{eq: Frobenius solutions at varphi=0}
    \varpi_j(\varphi)=\varphi\sum_{k}^3\binom{j}{k}F_k(\varphi)\log^{j-k}(\varphi)\, ,
\end{equation}
where $F_j(\varphi)$ are power series that are computed from the Picard-Fuchs equation. They are normalised so that $F_j(0)=\delta_{j,0}$.

We should stress that, when we work with the choice of integral symplectic basis and topological data in \eqref{eq: change of basis matrix at MUM point}, we are implicitly working on a free quotient of \eqref{eq: one parameter HV threefold}. As a result of this quotient, we will find that $\Gamma\in\mathbb{Q}^4$ in identities of the form \eqref{eq: generic integral of elliptic periods}.

\subsection{Normalisation of $\Pi$?}
Note that the periods in \eqref{eq: Frobenius solutions at varphi=0} have been normalised in such a way that 
\begin{equation}
    \varpi_0(\varphi)=\varphi+\mathcal{O}(\varphi^2).
\end{equation}
In other words, the Picard-Fuchs equation has indices $(1,1,1,1)$ at the large complex structure point $\varphi=0$. This normalisation has been included with the benefit of hindsight.

Recall that the holomorphic $3$-form at a given value of $\varphi$ is only defined up to multiplication by a non-zero constant. We are usually free to multiply the periods $\Pi$ by any non-vanishing holomorphic function of $\varphi$ without changing anything of consequence. However, a normalisation is implicitly determined by Equation \eqref{eq: generic integral of elliptic periods}, which we can find by analytically computing a single identity of the form \eqref{eq: generic integral of elliptic periods} as a function of $\varphi$. Alternatively, we may fix a number of values of $\varphi$, numerically compute a number of identities of the form \eqref{eq: generic integral of elliptic periods} and then guess what the normalisation of $\Pi$ must be. This is the quick and dirty approach that has been taken here.

\section{Identities}

We identify solutions of AESZ34 with identities of the form \eqref{eq: generic integral of elliptic periods}. We compute both $\Pi(\varphi)$ and the right hand side of \eqref{eq: generic integral of elliptic periods} and then identify the appropriate $\Gamma\in\mathbb{Q}^4$ by using the LLL-algorithm to look for a relation between the numerical values on the left and right hand sides of the identity. All identities in this section are found numerically at the fixed value of $\varphi=\frac{1}{64}$. However, we expect that they continue to hold as $\varphi$ is varied.

\subsection{Monodromy of elliptic periods}

We expect identities of the form \eqref{eq: generic integral of elliptic periods} when we can argue that the $3$-chain on the right-hand side is really a $3$-cycle. In order to argue that this is the case, we need the monodromy matrices of $\omega(\lambda,1)\otimes \omega(\lambda,\varphi)$ around the various singularities in the complex $\lambda$ plane.

We use \eqref{eq: integral symplectic basis of HV ellipitc curves} to fix a basis of elliptic periods at $\lambda=0$ and $\varphi=\frac{1}{64}$. We then analytically continue to the various singularities along the upper half $\lambda$-plane and compute anti-clockwise monodromy matrices around the various singularities. Around a singularity $\lambda_*$, the elliptic periods transform as
\begin{equation}
\omega(\lambda,1)\otimes\omega(\lambda,\varphi)\rightarrow \mu_{\lambda_*}~\omega(\lambda,1)\otimes\omega(\lambda,\varphi)\, .
\end{equation}
%\newpage
We find that
\begin{align}
    \mu_0 &= \begin{pmatrix} 1 & 0 \\ 3 & 1\end{pmatrix}\otimes \begin{pmatrix} 1 & 0 \\ 3 & 1\end{pmatrix}\\[4pt]
    \mu_{\frac{\varphi}{(1+2\sqrt{\varphi})^2}}&=I_2\otimes \begin{pmatrix}1 & -2\\ 0 & 1\end{pmatrix}\\[4pt]
    \mu_\varphi&=I_2\otimes\begin{pmatrix} 5 & -4\\ 4 & -3\end{pmatrix}\\[4pt]
    \mu_{\frac{\varphi}{(1-2\sqrt{\varphi})^2}}&=I_2\otimes \begin{pmatrix}5 & -2\\ 8 & -3\end{pmatrix}\\[4pt]
    \mu_{\frac{1}{9}}&=\begin{pmatrix} 1 & -2\\ 0 & 1\end{pmatrix}\otimes I_2\\[4pt]
    \mu_1 &=\begin{pmatrix} 7 & -6\\ 6 & -5\end{pmatrix}\otimes I_2~.
\end{align}

\subsection{Vanishing cycle at $\varphi=\frac{1}{9}$}
We expect that the (non-trivial) invariant $1$-cycle of a $2\times 2$ monodromy matrix is the vanishing cycle at the associated singularity. As with Hadamard products, we can use this to identify $3$-cycles in the fibred product of elliptic surfaces.\cite{Elmi:2023hof} For example
\begin{equation}
    \begin{pmatrix} 1 & -2\\ 0 & ~1 \end{pmatrix} \begin{pmatrix}1 \\ 0 \end{pmatrix}=\begin{pmatrix} 1 \\ 0 \end{pmatrix},
\end{equation}
which tells us that $(1,0)^T$ is a vanishing cycle of the first elliptic curve at $\lambda=\frac{1}{9}$.
 Similarly, the second elliptic curve has vanishing cycle $(1,1)^T$ at $\lambda=\varphi$. This leads us to the identity 
\begin{equation}
    \frac{1}{2}(1,0,-5,1)\Sigma_4\Pi(\varphi)=(1,0)\otimes(1,1)\Sigma_{2,2}\int_{\varphi}^\frac{1}{9}\omega(\lambda,1)\otimes\omega(\lambda,\varphi)d\lambda\, .
\end{equation}
As $\varphi\rightarrow \frac{1}{9}$, the integration contour in the above identity vanishes, so we expect that $(1,0,-5,1)^T$ is a vanishing cycle at $\varphi=\frac{1}{9}$. This is indeed the case. However, we must be careful. The precise vanishing cycle we compute depends on the branch that we approach the singularity on. Different choices of branches simply multiply the vanishing cycle by a monodromy matrix. In this case, the vanishing cycle computed in \cite{Candelas:2019llw} is $\left(1,0,5,1\right)^T$, which is simply
% \begin{equation}
%     \begin{pmatrix} 1 \\ 0 \\ 5 \\1 \end{pmatrix},
% \end{equation}
% which is simply
\begin{equation}
    \begin{pmatrix} 1 \\ 0 \\ 5 \\1 \end{pmatrix} = M_{\frac{1}{25}}\begin{pmatrix} \, 1 \\ \, 0 \\ -5 \\\, 1 \end{pmatrix}\, ,
\end{equation}
where $M_{\frac{1}{25}}$ is the monodromy matrix around $\varphi=\frac{1}{25}$ that was also computed in \cite{Candelas:2019llw}.

\subsection{Vanishing cycle at $\varphi=\frac{1}{25}$}
Similarly, we find that 
\begin{equation}
    \frac{1}{2}(0,0,-5,0)\Sigma_4\Pi(\varphi)=(1,0)\otimes(1,2)\Sigma_{2,2}\int_{\frac{\varphi}{(1-2\sqrt{\varphi})^2}}^{\frac{1}{9}} \omega(\lambda,1)\otimes\omega(\lambda,\varphi)d\lambda.
\end{equation}
Again, the contour of integration vanishes as $\varphi\rightarrow \frac{1}{25}$, so we would expect that $(0,0,-5,0)^T$ is a vanishing cycle at $\frac{1}{25}$. This is indeed the case. It is just a multiple of the vanishing cycle computed in \cite{Candelas:2019llw}.

\subsection{Holomorphic period at $\varphi=0$ and $T^3$}
As was the case with Hadamard products, we can identify non-trivial $3$-cycles by taking the union of some monodromy invariant $2$-cycle in the fibre over a contour on the base.\cite{Elmi:2023hof} In this way, we can identify a $3$-torus ($T^3$) that can be integrated over to obtain the holomorphic period at $\varphi=0$. More precisely, we find that
\begin{equation}
    (-1,0,0,0)\Sigma_4\Pi(\varphi)=(0,1)\otimes(0,1)\Sigma_{2,2}
    \oint_{\ell_0\,\ell_{\frac{\varphi}{(1+2\sqrt{\varphi})^2}}\,\ell_\varphi\,\ell_{\frac{\varphi}{(1-2\sqrt{\varphi})^2}}}\omega(\lambda,1)\otimes\omega(\lambda,\varphi)d\lambda\, ,
\end{equation}
where the integration contour 
\begin{equation*}
\ell_0\,
\ell_{\frac{\varphi}{(1+2\sqrt{\varphi})^2}}\,
\ell_\varphi\,
\ell_{\frac{\varphi}{(1-2\sqrt{\varphi})^2}}
\end{equation*}
is a closed contour on the $\lambda$ plane that first encircles $\lambda=0$ in an anti-clockwise fashion, then encircles $\lambda=\frac{\varphi}{(1+2\sqrt{\varphi})^2}$
in an anti-clockwise fashion, then encircles... . We see from the above monodromy matrices that $(0,1)\otimes(0,1)$ is invariant under monodromy around this contour. We have therefore identified a $T^3$ we can integrate over and obtain the holomorphic period at the point of maximal unipotent monodromy $\varphi=0$. The results of \cite{Elmi:2023hof} suggest that this family of $T^3$s defines an SYZ-like fibration by constant phase contours on the $\lambda$-plane.

\bibliographystyle{utcaps}
\bibliography{bibliography}
 
\end{document}